\documentstyle[amssymb]{article}
\date{}
\begin{document}
{\noindent\Large
\texttt{On the $p$-th root of a $p$-adic number}
\footnote{AMS MCS 11S99, 13K05.}}

\bigskip\noindent
{{Alfonso Di Bartolo}\footnote{Supported by Universit\`a di Palermo (Co.R.I.).}\\
\footnotesize{Dipartimento di Matematica e Applicazioni}\vspace*{-1,5mm}\\
\footnotesize{Via Archirafi 34, I-90123 Palermo (Italy)}\vspace*{-1,5mm}\\
\footnotesize{alfonso@math.unipa.it}}

\medskip\noindent
{{Giovanni Falcone}\footnote{Supported by M.I.U.R., Universit\`a di Palermo  (Co.R.I.).}\\
\footnotesize{Dipartimento di Metodi e Modelli Matematici}\vspace*{-1,5mm}\\
\footnotesize{Viale delle Scienze Ed. 8, I-90128 Palermo (Italy)}\vspace*{-1,5mm}\\
\footnotesize{gfalcone@unipa.it}}

\vskip1cm\begin{abstract} \noindent We give a
sufficient and necessary condition for a $p$--adic integer to have
$p$--th root in the ring of $p$--adic integers. The same condition
holds clearly for residues modulo $p^k$. We give a proof that
Fermat's last theorem is false for $p$--adic integers and for
residues mod $p^k$.
\end{abstract}

\large\vskip1cm\noindent 
Under the assumption that the prime $p$ does
not divide the integer $k$, an immediate consequence of Hensel's
lemma is that a $p$--adic unit $a=l_0+pl_1+p^2l_2+\cdots $ has a
$k$--th root in the ring of $p$--adic integers if and only if
$l_0$ has a $k$--th root in ${\Bbb
Z}_p={\mathbb{Z}}/p{\mathbb{Z}}$. The same argument gives for the
$p$--th root a sufficient but not necessary condition. In order to
find the $p$--th root, we apply the exponential and logarithm maps
to the Witt ring ${\mathfrak W}({\mathbb Z}_p)$, which is
isomorphic to the ring of $p$--adic integers.

In his paper \cite{Witt} of 1936, E. Witt found the algorithm
which gives recursively the factor systems necessary to describe the ring
of $p$--adic integers as the inverse limit of the rings of residues ${\mathbb{Z}}_{p^k}={\mathbb{Z}}/p^k{\mathbb{Z}}$.
In this context this ring is denoted by
$${\mathfrak W}({\mathbb Z}_p)=\{{\mathbf{x}}=(x_0,x_1,\cdots,x_k,\cdots)
 |\, x_i\in{\mathbb{Z}}_p\}$$
and its elements are called \em Witt vectors\em . For a detailed
exposition, we refer to \cite{DMG}, Ch.~V, no.~1.

The ground subring generated by the unitary element ${\mathbf
1}=(1,0,0,\cdots)$ is isomorphic to $\mathbb Z$ and in \cite{JNT}
we gave the representation of an arbitrary natural integer
$n\in{\mathbb N}$ as the element $n{\mathbf 1}$ of ${\mathfrak
W}({\mathbb Z}_p)$.

One of the reasons to use the representation of
integers as Witt vectors is that the quotient ring ${\mathfrak W}_k({\mathbb Z}_p)={\mathfrak W}({\mathbb Z}_p)/p^k{\mathfrak W}({\mathbb Z}_p)$,
which is isomorphic to the ring ${\mathbb{Z}}/p^k{\mathbb{Z}}$ of
residues modulo $p^k$, can be represented as the ring of truncations
of Witt vectors after the first $k$ entries, that is the set of
elements of the shape
$$(x_0,x_1,\cdots ,x_{k-1}].$$
This allows one to consider simultaneously integers, rationals,
$p$--adics and residues modulo $p^k$ in many arguments.

\bigskip\noindent{\textbf{1. Integers, residues and $p$--adics as Witt
vectors.}}

\smallskip\noindent
For any $k=0,1,\cdots$, let $a_0,\cdots ,a_k\in {\mathbb
Q}$ be such that
$$\Phi_k(a_0,\cdots ,a_k)=a_0^{p^k}+pa_1^{p^{k-1}}
+\cdots +p^ka_k= n.$$ Then we have (cf. \cite{JNT}): \begin{itemize} \item[i)] $a_0=n$
and
$a_{k+1}=\sum_{i=0}^k{1\over{p^{k-i+1}}}(a_i^{p^{k-i}}-a_i^{p^{k-i+1}})$ are integers;
\item[ii)] $n\cdot {\mathbf 1} =
({\overline{a}}_0,{\overline{a}}_1,\cdots )$;
\item[iii)] if $p$ does
not divide $n$, then $n$ divides each $a_k$.\end{itemize}
 We identify any integer $n$, not divisible by $p$,
with
$$n\equiv n\cdot {\mathbf{1}}=(n,-n{\sf q}_1(n),-n{\sf q}_2(n),\cdots ),$$
where we take the entries modulo $p$ and we put ${\sf
q}_i(n)=-a_i/n$, which is an integer by \em iii\em ) in the above
Proposition. We remark that ${\sf q}_1(n)={\frac{n^{p-1}-1}{p}}$ is
the \em Fermat quotient \em of $n$.

Furthermore we identify the residue of $n$ modulo $p^k$
with the truncated Witt vector $(n,-n{\sf q}_1(n),-n{\sf q}_2(n),\cdots ,-n{\sf q}_{k-1}(n)]$.

Let now $$a=\sum_{i=0}^\infty l_ip^i$$ be a $p$--adic unit, with $0< l_0 <p$ and
$0\leq l_i <p$, for $i>0$. The first two entries of the Witt vector
corresponding to $a$ are therefore the same of $l_0+l_1\cdot p$,
that is
$$a\equiv (l_0,-l_0{\sf
q}_1(l_0),\cdots)+(0,l_1,\cdots )+\cdots =(l_0,l_1-l_0{\sf
q}_1(l_0),\cdots ).$$ In the next paragraph we will
compute the first entries of the Witt vector
corresponding to a rational number.\hfill$\Box$

\medskip\noindent According
to a consolidated notation (\cite{DMG}, Ch.~5, {\S\S} 1, 3, 4), the
element $(\bar{a},0,0,0,\cdots )\in {\mathfrak W}({\mathbb Z}_p)$ is denoted by
${a}^\tau $ and is called the \em Teichm\"uller representative of $a$\em . Any Witt vector
${\mathbf{x}}=(x_0,x_1,x_2,\cdots )$ such that $x_0\not\equiv 0\;(\mbox{mod }p)$
can be written as the product ${\mathbf{x}}=x_0^\tau(1,x_1/x_0,x_2/x_0,\cdots )$. The invertible elements in ${\mathfrak W}({\mathbb Z}_p)$
are precisely those ${\mathbf{x}}=(x_0,x_1,x_2,\cdots )$ having $x_0\not\equiv 0\;(\mbox{mod }p)$. Therefore any element of
the quotient field of ${\mathfrak W}({\mathbb Z}_p)$, which is (isomorphic to) the field of $p$--adic numbers,
can be written as ${\mathbf{x}}=p^zx_0^\tau(1,x_1/x_0,x_2/x_0,\cdots )$,
with $z\in {\mathbb{Z}}$ and $x_0\not= 0$. The rational integer $p^{-z}$ is the $p$--adic valuation $|{\mathbf{x}}|_p$ of ${\mathbf{x}}$. With a slight abuse, we will call such elements \em Witt vectors\em , as well.\hfill$\Box$

\bigskip\noindent{\textbf{2. Logarithm and exponential map. De Moivre
formula.}}

\smallskip\noindent
In this paragraph we assume $p>2$. The formal power series
\begin{itemize}
  \item[] $\log(1+p{\textbf{x}})=p{\textbf{x}}-1/2(p{\textbf{x}})^2+1/3(p{\textbf{x}})^3
-\cdots $
  \item[] $\mbox{e}^{p{\textbf{x}}}=1+p{\textbf{x}}+1/2!(p{\textbf{x}})^2+1/3!(p{\textbf{x}})^3
+\cdots $
\end{itemize}
are simply polynomials in the ring ${{\mathfrak W}_k}({\mathbb
Z}_p)$ of truncated Witt vectors, isomorphic to the ring of
residues ${\mathbb{Z}}_{p^k}={\mathbb{Z}}/p^k{\mathbb{Z}}$. For
instance, for $p>3$, we have
\begin{itemize}
  \item[] $\log(1,a_1,a_2]=(0,a_1,a_2-{\frac{1}{2}}a_1^2]$,
  \item[] $\mbox{e}^{(0,a_1,a_2]}=(1,a_1,a_2+{\frac{1}{2}}a_1^2]$.
\end{itemize}
Since the two maps can be defined for any $k>0$, they are defined
on the whole of
\begin{itemize}
  \item[] ${\mathbf{1}}+p\, {\mathfrak W}({\mathbb Z}_p)=\{{\mathbf{x}}=(1,x_1,x_2,\cdots ):x_i\in{\Bbb{Z}}_p\}$ and
  \item[] $p\, {\mathfrak W}({\mathbb Z}_p)=\{{\mathbf{x}}=(0,x_1,x_2,\cdots ):x_i\in{\Bbb{Z}}_p\}$,
\end{itemize}
respectively, and the two maps are mutually inverse.

\smallskip
Let ${\mathbf{x}}=p^zx_0^\tau(1,x_1/x_0,x_2/x_0,\cdots )$,
with $z\in {\mathbb{Z}}$ and $x_0\not\equiv 0\;(\mbox{mod }p)$, be an arbitrary
Witt vector. If we define
\begin{itemize}
  \item[] the module $\rho_{\mathbf{x}}:=p^zx_0^\tau$,
  \item[] the argument $\vartheta_{\mathbf{x}}:=
  \log(1,x_1/x_0,x_2/x_0,\cdots )\in {\mathfrak{W}}$,
\end{itemize}
then we can write $${\mathbf{x}}=
\rho_{\mathbf{x}}\mbox{e}^{\vartheta_{\mathbf{x}}}$$ and recover De Moivre
formula
\begin{itemize}
  \item[] $\rho_{\mathbf{xy}}=\rho_{\mathbf{x}}\rho_{\mathbf{y}}$,
  \item[] $\vartheta_{\mathbf{xy}}=
  \vartheta_{\mathbf{x}}+\vartheta_{\mathbf{y}}$,
\end{itemize}
holding for $p$--adics as well as for residues modulo
$p^k$. We remark that, modulo $p^2$, De Moivre formula $\vartheta_{nm}=
  \vartheta_{n}+\vartheta_{m}$ coincides with the Eisenstein congruence
${\sf q}_1(n\cdot m)\equiv {\sf q}_1(n)+{\sf q}_1(m)\;(\mbox{mod
}p)$.\hfill$\Box$

\medskip\noindent As an application we compute
$${\mathbf{x}}^{-1}=\rho_{\mathbf{x}}^{-1}\mbox{e}^{-\vartheta_{\mathbf{x}}}$$
for a natural integer $n\equiv n^\tau (1,-{\sf q}_1(n),-{\sf
q}_2(n),\cdots )$, not divisible by $p$. In fact,
$$\Big(n^\tau (1,-{\sf q}_1(n),-{\sf q}_2(n),\cdots )\Big)^{-1}=
(n^\tau)^{-1}\mbox{e}^{-(0,-{\sf q}_1(n),-{\sf q}_2(n)-{\frac{1}{2}}{\sf
q}_1^2(n),\cdots )}=$$
$$(n^{-1})^\tau \mbox{e}^{(0,{\sf q}_1(n),{\sf q}_2(n)+{\frac{1}{2}}{\sf
q}_1^2(n),\cdots )}=(n^{-1})^\tau {(1,{\sf q}_1(n),{\sf q}_2(n)+{\sf
q}_1^2(n),\cdots )}.$$ Similarly, if $m$ and $n$ are two integers, not divisible by $p$, we find
$${\frac{m}{n}}\equiv\left({\frac{m}{n}},-{\frac{m}{n}}({\sf
q}_1(m)-{\sf q}_1(n)),\cdots \right).$$
 \hfill$\Box$

\medskip\noindent It is standard to define, for ${\mathbf{x}}\in 1+p\, {\mathfrak W}({\mathbb Z}_p)$ and ${\mathbf{y}}\in{\mathfrak W}({\mathbb Z}_p)$,
$${\mathbf{x}}^{\mathbf{y}}:=\mbox{e}^{{\mathbf{y}}\log{\mathbf{x}}}\in{\mathfrak W}({\mathbb Z}_p),$$
and the aim of this paper is to remark that ${\mathbf{x}}^{\mathbf{y}}$ is still in ${\mathfrak W}({\mathbb Z}_p)$ for a $p$--adic number ${\mathbf{y}}$ with positive $p$--adic valuation $|{\mathbf{y}}|_p= p^{k}$, if we assume $x_i\equiv 0\;(\mbox{mod }p)$ for $i=1,2,\cdots ,k$.\hfill$\Box$

\bigskip\noindent{\textbf{3. The $p$--th root.}}

\medskip\noindent
Let $p>2$ and let ${\mathbf{x}}=p^z
x_0^\tau(1,x_1/x_0,\cdots,x_k/x_0,\cdots)$ be a Witt vector, with
$z\in {\mathbb{Z}}$ and $x_0\not\equiv 0\;(\mbox{mod }p)$. As an
immediate consequence of De Moivre formula, we have
$${\mathbf{x}}^{p^k}= p^{zp^k} x_0^\tau(1,\underbrace{0,\cdots ,0}_{k},x_1/x_0,\cdots)$$
(note that, from the $k+2$--nd one on, the entries become more involved).
Furthermore, if the Witt vector
${\mathbf{x}}=(x_0,x_1,\cdots)$ is such that $x_0\not\equiv 0$ and $x_i\equiv 0$, for $i=1,\cdots ,k $,
then we find
$${\frac{{\mathbf{x}}-{\mathbf{x}}^p}{p^{k+1}}}\equiv{\frac{(x_0,0,\cdots,0,x_{k+1}]-(x_0,0,\cdots,0,0]}{p^{k+1}}}$$
$$={\frac{(0,0,\cdots,0,x_{k+1}]}{p^{k+1}}}=(x_{k+1},\cdots],$$
(once again, we remark that, from the $k+2$--nd one on, the entries become more involved). Thus we have in this case $$x_{k+1}\equiv -{\frac{1}{p^{k}}}{\frac{{\mathbf{x}}^p-{\mathbf{x}}}{p}}\;(\mbox{mod }p)=
-{\frac{1}{p^{k}}}{\mathbf{x}}\, {\sf q}_1({\mathbf{x}})\;(\mbox{mod }p),$$
where, in analogy to the case of an integer, we define the \em Fermat quotient \em of the Witt vector
 ${\mathbf{x}}=(x_0,x_1,\cdots,x_k,\cdots)$, having $x_0\not\equiv 0\;(\mbox{mod }p)$, as ${\sf q}_1({\mathbf{x}})
 ={\frac{{\mathbf{x}}^{p-1}-1}{p}}\in {\mathfrak W}({\mathbb Z}_p)$. We note that, in accordance with the case of an integer, we have $x_1\equiv -{\mathbf{x}}\,{\sf q}_1({\mathbf{x}})\;(\mbox{mod }p)$ and again, De Moivre formula $\vartheta_{\mathbf{xy}}=
  \vartheta_{\mathbf{x}}+\vartheta_{\mathbf{y}}$ reduces in ${\mathfrak W}_2({\mathbb Z}_p)$ to Eisenstein congruence
  ${\sf q}_1({\mathbf{x}}\cdot{\mathbf{y}})\equiv {\sf q}_1({\mathbf{x}})+{\sf q}_1({\mathbf{y}})\;(\mbox{mod }p)$.\hfill$\Box$

\medskip\noindent
Having the entries $x_i\equiv 0\;(\mbox{mod }p)$ for $i=1,\cdots ,
k$ is not only a necessary condition for a Witt vector to be a
$p^k$--th power, it is sufficient, as well. Our condition is based
on the fact that
$$p^{-k}\log(1,x_1/x_0,\cdots,x_k/x_0,\cdots)=p^{-k}(0,x_1/x_0,\cdots)$$
lies in $p\, {\mathfrak W}({\mathbb Z}_p)$ if and only if
$x_i\equiv 0\;(\mbox{mod }p)$, for $i=1,2,\cdots k$. The Witt
vector ${\mathbf{x}}=p^z
x_0^\tau(1,x_1/x_0,\cdots,x_k/x_0,\cdots)$ has therefore a $p^k$--th root in
${\mathfrak W}({\mathbb Z}_p)$ if and only if $z\equiv
0\;(\mbox{mod }p^k)$ and $x_i\equiv 0\;(\mbox{mod }p)$, for
$i=1,2,\cdots k$. In this case the root is unique and it is
$${\mathbf{x}}^{{\frac{1}{p^k}}}=p^{{\frac{z}{p^k}}}x_0^\tau \mbox{e}^{{\frac{1}{p^k}}\log(1,0,\cdots ,0,x_{k+1}/x_0,\cdots
)}.$$
For instance, let ${\mathbf{x}}=x_0^\tau(1,0,x_2/x_0,\cdots,x_k/x_0,\cdots)\in{\mathfrak W}({\mathbb Z}_p)$. Then we have
$${\mathbf{x}}^{{\frac{1}{p}}}\equiv x_0^\tau \mbox{e}^{{\frac{1}{p}}\log(1,0,x_2/x_0,x_3/x_0]}=x_0^\tau \mbox{e}^{(0,x_2/x_0,x_3/x_0]}$$
$$=x_0^\tau {(1,x_2/x_0,x_3/x_0+1/2(x_2/x_0)^2]}\;(\mbox{mod }p^3).$$
Therefore the integer $n$, not divisible by $p$, has the $p$--th root in the ring of $p$--adic integers if and only if ${\sf q}_1(n)\equiv 0
\;(\mbox{mod }p)$, that is if $n^p\equiv n\;(\mbox{mod
}p^2)$ and the $p$--adic unit $a=l_0+pl_1+p^2l_2+\cdots $ has the $p$--th root in the ring of $p$--adic integers if and only if
${\sf q}_1(a)\equiv 0
\;(\mbox{mod }p)$, that is if
$a^p\equiv a\;(\mbox{mod }p^2)$. We remark that this condition is equivalent to say that $l_1\equiv {\frac{l_0^p-l_0}{p}}\;(\mbox{mod }p)$.
\hfill$\Box$

\medskip\noindent
If $p=2$, the two opposite square roots of a unit ${\mathbf{x}}$
exist if and only if ${\mathbf{x}}\equiv 1\;(\mbox{mod }8)$. This
follows directly from Hensel's lemma. But we note that it is
possible to compute these roots also as
${\mathbf{x}}^{{\frac{1}{2}}}=
\mbox{e}^{{\frac{1}{2}}\log{\mathbf{x}}}$. In fact, it is
well--known that for $p=2$ the exponential map is defined for
${\mathbf{x}}\in 4\, {\mathfrak W}({\mathbb Z}_2)$.\hfill$\Box$

\medskip\noindent
Let ${\mathbf{x}}=(x_0,x_1,\cdots,x_k,\cdots)$ be a Witt vector, such that $x_0\not\equiv 0$ and $x_i\equiv 0\;(\mbox{mod }p)$, for $i=1,2,\cdots ,k$, and compute $${\mathbf{x}}^{{\frac{1}{p^k}}}\equiv (x_0,x_{k+1}]\;(\mbox{mod }p^2).$$
Thus the above congruence $x_{k+1}\equiv
-{\frac{1}{p^{k}}}{\mathbf{x}}\, {\sf q}_1({\mathbf{x}})\;(\mbox{mod }p)$ can be written meaningfully as
$${\sf q}_1\Big({\mathbf{x}}^{{\frac{1}{p^k}}}\Big)\equiv{\frac{1}{p^{k}}}\, {\sf q}_1({\mathbf{x}})\;(\mbox{mod }p),$$
in accordance with the Eisenstein congruence ${\sf q}_1({\mathbf{x}}\cdot {\mathbf{y}})\equiv {\sf q}_1({\mathbf{x}})+{\sf q}_1({\mathbf{y}})\;(\mbox{mod }p)$.\hfill$\Box$

\medskip\noindent \emph{Example}: A non--trivial case where ${\sf q}_1(n)\equiv 0\;(\mbox{mod }p)$
is for $n=3$ and $p=11$. This means that $n=3$ has $11$--adic root in the $11$--adic field or, equivalently, that
the residue of $n=3$ in ${\mathbb{Z}}_{11^k}$ has
$11$--th root in ${\mathbb{Z}}_{11^k}$, for any $k\geq 1$. In particular, we find
$$3^{{\frac{1}{11}}}=3^\tau \mbox{e}^{({\frac{1}{11}}\log(1,0,-{\sf q}_2(3)])}=3^\tau \mbox{e}^{({\frac{1}{11}}(0,0,-{\sf q}_2(3)])}$$
$$=3^\tau \mbox{e}^{(0,-{\sf q}_2(3)]}=3^\tau {(1,-{\sf q}_2(3)]}.$$
As we mentioned above, a consequence of the congruence ${\sf q}_1(3)\equiv 0\;(\mbox{mod }11)$
is that ${\sf q}_2(3)\equiv \frac{{\sf q}_1(3)}{11}\equiv
5368\equiv 4\;(\mbox{mod }11)$. Therefore
$$3^{{\frac{1}{11}}}\equiv 3^\tau {(1,-4]}=(3,-1]$$
hence $3-11=-8$ is the $11$--th root of $3$ modulo
$11^2$.\hfill$\Box$

\medskip\noindent Denote by $\varphi_1(x_0,y_0)$ the factor system
defining the sum in the ring ${\mathfrak W}_2({\mathbb Z}_p)$ of truncated Witt vectors, that
is
$$(x_0,x_1]+(y_0,y_1]=(x_0+y_0,x_1+y_1+\varphi_1(x_0,y_0)].$$
As remarked in \cite{JNT}, we have
$$\varphi_1(x_0,y_0)\equiv\sum_{i=1}^{p-1}{\frac{(-1)^i}{i}}x_0^iy_0^{p-i}\;(\mbox{mod }p).$$
The smallest prime $p$ such that, for a suitable integer $0<x<p-1$,
$$\varphi_1(1,x)\equiv 0\;(\mbox{mod }p)$$ is $p=7$. In fact,
$\varphi_1(1,2)\equiv 0\;(\mbox{mod }7)$. Since
$$129=1^7+2^7\equiv(1,0]+(2,0]=(3,0]\;(\mbox{mod }7^2),$$ it follows that
$129$ is the $7$--th power of a $7$--adic integer. This shows that
the equality $x^7+y^7+z^7=0$ has a non--trivial $7$--adic solution
and the equality $x^7+y^7+z^7\equiv 0\;(\mbox{mod }7^k)$ has a
non--trivial solution for any $k\geq 0$ (cfr. \cite{BoSch}, Remark
1, p. 163).\hfill$\Box$

\medskip\noindent It seems very rare that $n=2$ has the $p$--th root in the field of
$p$--adics. In fact $1093$ and $3511$ are the only known primes, up
to $1.25\, \cdot 10^{15}$, for which ${\sf q}_1(2)\equiv 0\;(p)$.
These primes are called \em Wieferich primes \em since Wieferich
proved in $1909$ that, if $x^p+y^p+z^p=0$ had a non trivial integer solution
with $xyz$ not divisible by $p$, then ${\sf q}_1(2)\equiv 0\;(p)$. In $1910$ Mirimanoff proved moreover that for such a prime $p$ it must hold
that ${\sf q}_1(3)\equiv 0\;(p)$ and a still open question is
whether it is possible that simultaneously ${\sf q}_1(2)\equiv {\sf
q}_1(3)\equiv 0\;(p)$.\hfill$\Box$

\end{document}